\documentclass[10pt]{article}

\usepackage{amsmath}

\newtheorem{theorem}{Theorem}

\title
{
    On a Problem Posed by Maurice Nivat
}

\author
{
    Maxim A. Babenko
    \thanks
    {
		Dept. of Mechanics and Mathematics,
		Moscow State University, Vorob'yovy Gory, 119899 Moscow,
		Russia, \textsl{email}: mab@shade.msu.ru.
    }
}

\begin{document}

\maketitle

\begin{abstract}

Consider a $m \times n$ matrix~$A$, whose elements are arbitrary
integers. Consider, for each square window of size $2 \times
2$, the sum of the corresponding elements of $A$. These sums form
a $(m - 1) \times (n-1)$ matrix~$S$. Can we
efficiently (in polynomial time) restore the
original matrix~$A$ given~$S$?

This problem was originally posed by Maurice Nivat
for the case when the elements of matrix~$A$ are zeros and ones.
We prove that this problem is solvable in polynomial time.
Moreover, the problem still can be efficiently solved if
the elements of $A$ are integers from given intervals.

On the other hand, for $2 \times 3$ windows the
similar problem turns out to be NP-complete.

\end{abstract}

\section{Introduction}

By $M_{mn}$ we denote the set of all $n \times m$ integer matrices.
For a given matrix $A \in M_{mn}$ consider the matrix of sums
for all $2 \times 2$ windows, denoted by~$S = \Sigma_{22}(A)$.
Here indices denote the window size. More generally,
a mapping $\Sigma_{m' n'} : M_{mn} \to M_{m-m'+1,n-n'+1}$ (for $m'\times n'$ window)
is defined in a similar way.

Now let $A$ be a matrix with 0-1 elements and $S =
\Sigma_{22}(A)$. How can we restore $A$ knowing $S$? First of all,
note that $S$ could have many preimages. (For example, consider an
arbitrary 0-1 matrix such that every its column is formed by
alternating zeros and ones. Clearly all elements of $S$ are equal
to 2.) So our goal is to find (efficiently) one of the preimages of
$S$ if they exist.

We also consider a more general problem with \emph{upper constraints}.
Namely, given a matrix $S$ of $m' \times n'$
sums and \emph{upper constraints matrix}~$U \in M_{mn}$ we look for a matrix~$A$
such that
\begin{eqnarray*}
    0 \le A \le U, \\
    \Sigma_{m'n'}(A) = S.
\end{eqnarray*}
As usual, $0 \le A \le U$ means that $0 \le A_{ij} \le U_{ij}$ for all $i$ and $j$.

The original problem (with $U_{ij} = 1$) is called \emph{binary}.
In this paper we prove the following results:

\begin{theorem}
\label{binary_2x2}
    The binary problem with $2 \times 2$ window is solvable in polynomial time.
    Moreover, in a typical RAM model it can be solved in $O(mn)$ time.
\end{theorem}

\begin{theorem}
\label{upper_2x2}
    In a typical RAM model the problem with $2 \times 2$ window and upper constraints
    is solvable in $O(mn (m + n) (1 + U_{min}))$ time, where $U_{min} = \min_{ij} U_{ij}$.
\end{theorem}

Suppose the elements of $U$ and $S$ are given in unary notation.
Then Theorem~\ref{upper_2x2} implies that the binary
problem with upper constraints is solvable in polynomial time.
The next theorem shows the hardness of the similar problem for
$2 \times 3$ window.

\begin{theorem}
\label{upper_2x3}
    The problem with upper constraints (given in unary notation)
    and $2 \times 3$ window is NP-complete.
\end{theorem}

\section{Binary Problem for $2 \times 2$ Window}

Let $A$ be a matrix we are looking for and $S$ be the matrix of
sums that is given to us. We number the rows and the columns
starting from zero (rows~$0, \dots, m - 1$ and columns~$0, \dots,
n - 1$).

Note, that it is sufficient to restore only the elements in the
zero row and column of $A$. After that, all other elements are
determined uniquely. We start with an observation that works not
only in the binary case ($A_{ij} \in \{ 0, 1 \}$), but also for
any upper constraints ($0 \le A_{ij} \le U_{ij}$).

We may assume that $A_{00}$ is already known (since we
can consider all $U_{00} + 1$ possible cases one by one). Let
$x_1,\ldots,x_{n-1}$ be the remaining elements of the zero row of $A$ and
$y_1,\ldots,y_{m-1}$ be the remaining elements of the zero column:
\begin{center}
    \begin{tabular}{|c|c|c|c|c|}
        \hline
        $A_{00}$  & $x_1$ & $x_2$ & \ldots & $x_{n-1}$ \\
        \hline
        $y_1$     &       &       &        &           \\
        \hline
        $y_2$     &       &       &        &           \\
        \hline
        \vdots    &       &       &        &           \\
        \hline
        $y_{m-1}$ &       &       &        &           \\
        \hline
    \end{tabular}
\end{center}

Easy induction shows that
$$
    A_{ij} = (-1)^i x_j + (-1)^j y_i + b_{ij},
$$
where $b_{ij}$ are some constants depending on $A_{00}$ and matrix~$S$.
The numbers $b_{ij}$
can be computed in $O(mn)$ time.
So we get the following requirements for $x_j$ and $y_i$:
\begin{equation}
\label{ineq2x2}
    \begin{array}{lcccl}
        0 & \le & x_j & \le & U_{0j}; \\
        0 & \le & y_i & \le & U_{i0}; \\
        -b_{ij} & \le & (-1)^i x_j + (-1)^j y_i & \le & U_{ij} - b_{ij}.
    \end{array}
\end{equation}

Moreover, if conditions~(\ref{ineq2x2}) are satisfied for some
$x_j$ and $y_i$, then corresponding matrix~$A$ provides a solution
for the original problem with upper constraints.

Our algorithm uses that each inequality in (\ref{ineq2x2}) depends
on at most two variables. Suppose we consider the binary case.
Then $x_j$ and $y_i$ are Boolean variables and the
inequalities~(\ref{ineq2x2}) can be written as a Boolean formula.
Indeed, for each pair~$(i,j)$ the inequality $-b_{ij} \le (-1)^i
x_j + (-1)^j y_i \le 1 -b_{ij}$ forbids some pairs of
values~$(x_j, y_i)$. Putting these restrictions together we obtain
a 2-CNF formula in $x_j$, $y_i$. It is clear that the size of this
formula is $O(mn)$.

A well-known fact is that for a given 2-CNF formula one can find
whether it is satisfiable or not in polynomial time (and find a
satisfying assignment if it exists). This problem is often called
\emph{2-SAT problem}. Moreover, there exists an algorithm solving
2-SAT that runs in linear time (in the length of the formula). Our
formula is of $O(mn)$ size, and hence we obtain the proof of
Theorem~\ref{binary_2x2}.

In the rest of the section we briefly outline the idea behind the
linear time algorithm for solving 2-SAT problem. Let $\{z_i\}$ be
the set of Boolean variables. A \emph{literal} is a variable~$z_i$
(denoted by $z_i^0$) or its negation (denoted by $z_i^1$).
By \emph{2-CNF} we mean a formula~$\phi$ in conjunctive normal form where
each clause is a disjunction of at most two literals. Without loss of
generality we may assume that every clause has exactly two
literals (maybe identical).

Converting the disjunctions into implications we get:
$$
    \phi(z_1, \ldots, z_n) = \bigwedge_i (
        z_{\alpha_i}^{\sigma_i} \to z_{\beta_i}^{\pi_i}
    ),
$$
where~$\sigma_i, \pi_i \in \{ 0, 1 \}$.

Our first step is to construct a directed graph~$G = \langle V, E
\rangle$, where $V$ is the set of literals:
$$
    V = \{ z_i^0, z_i^1 \}.
$$
For each implication~$u \to v$ (where $u, v$ are literals) we add
arcs $u \to v$ and $\bar v \to \bar u$ (here
$\overline{z_i^\sigma}$ denotes $z_i^{1 - \sigma})$.

To satisfy~$\phi$ means to label vertices in this graph by Boolean values in
such a way that $z_i^0$ and $z_i^1$ get opposite values and there
is no arc going from a \textsc{true} vertex to a \textsc{false}
one.

The size of the graph is linear in the length of $\phi$.
We calculate the strongly-connected components of $G$
using depth-first search twice (see~\cite{CorLeiRiv1999}).
This requires linear time.

Suppose literals~$z_i^0$ and $z_i^1$ (for some~$i$) belong to the
same strongly-connected component. Then $\phi$ is unsatisfiable
since it implies both $z_i \to \bar{z}_i$ and $\bar{z}_i \to z_i$.

On the other hand, if literals~$z_i^0$ and
$z_i^1$ are in different components for each~$i$, then formula is
satisfiable. To show this we perform a topological sort of the components.
In other words, we assign natural numbers to the
components in such a way that for each arc going from
component~$C_i$ to component~$C_j$ we have $i \le j$.

Now we describe how to assign Boolean values to variables~$z_i$.
Consider a pair of literals $z_k$ and $\bar{z}_k$. Let $C_i$ be
the component containing $z_k$ and $C_j$ be the component
containing $\bar{z}_k$. If $i < j$ then we put~$z_k =
\textsc{false}$. Otherwise $i > j$ since $z_k$ and $\bar{z}_k$ are
in different components. In this case put~$z_k = \textsc{true}$.
It remains to prove that these values satisfy the formula~$\phi$, i.e.,
that no arc goes from \textsc{true}
to \textsc{false}.

Suppose the contrary and let~$u \to v$ be such an arc (here
$u$ and $v$ are literals). Let~$C_i$ denote the component
containing $u$ and let $C_j$ be the component containing~$v$. Then
$i \le j$. Consider vertex~$\bar u$ and vertex $\bar v$. Let
$C_{i'}$ and $C_{j'}$ be their components. Since~$u =
\textsc{true}$ and $v = \textsc{false}$ we have $i' < i$ and $j' >
j$, hence $i' < j'$. On the other hand the graph contains the
arc~$\bar v \to \bar u$ that violates topological order. The
correctness of the algorithm is now established.

It is clear that using appropriate data structures this algorithm can be implemented in linear time.

\section{The Case of $2 \times 2$ Window and Arbitrary Upper Constraints}

Now suppose that $A_{ij}$ are integers in the range $0 \ldots
U_{ij}$. We use the fact that the problem can be reduced to
the set of inequalities~(\ref{ineq2x2}). As above, we consider
each all possibilities for $A_{00}$ separately.

We let~$x_j = (-1)^j \alpha_j$, $y_i = (-1)^{i+1} \beta_i$.
Then the inequalities~(\ref{ineq2x2}) become two-sided constraints
on $\alpha_j$, $\beta_i$ and the differences $\alpha_j - \beta_i$:
\begin{equation}
\label{ineq2x2alt}
    \begin{array}{lcccl}
        L_j^1    & \le & \alpha_j & \le & U_j^1, \\
        L_i^2    & \le & \beta_i & \le & U_i^2, \\
        L_{ij}^3 & \le & \alpha_j - \beta_i & \le & U_{ij}^3
    \end{array}
\end{equation}
for some $L_j^1$, $U_j^1$, $L_i^2$, $U_i^2$, $L_{ij}^3$, $U_{ij}^3$.
Consider a more general (and more ``uniform'') set of inequalities:
\begin{equation}
\label{ineq2x2uniform}
    \begin{array}{lcccl}
        L_j^1    & \le & \alpha_j - \theta & \le & U_j^1, \\
        L_i^2    & \le & \beta_i - \theta & \le & U_i^2, \\
        L_{ij}^3 & \le & \alpha_j - \beta_i & \le & U_{ij}^3.
    \end{array}
\end{equation}
These two systems of inequalities are either both consistent or both
inconsistent. Indeed, every integer solution~$(\alpha_j, \beta_i)$ of~(\ref{ineq2x2alt})
can be transformed into a solution of~(\ref{ineq2x2uniform})
by setting~$\theta = 0$. And visa versa, if $(\alpha_j, \beta_i, \theta)$
is an integer solution of~(\ref{ineq2x2uniform}),
then $(\alpha_j - \theta, \beta_i - \theta)$ is
an integer solution of~(\ref{ineq2x2alt}). Thus it is enough
to consider inequalities~(\ref{ineq2x2uniform}) only.

This set of inequalities has a form of \emph{difference
constraints}. Using Ford--Bellman algorithm
(see~\cite{CorLeiRiv1999}) we may find an integer solution
for~(\ref{ineq2x2uniform}) or establish that it does not exist in
$O(mn (m + n))$ time.

Namely, suppose we have a set of variables~$\{ z_i \}$ and a set
of difference constraints $z_i - z_j \le w_{ij}$ for some~$i, j$
and integer constants~$w_{ij}$. Our task is to find an integer
solution (if it exists) for this set of inequalities. To do so, we
consider a directed graph~$G = \langle V, E \rangle$ constructed
in the following way. Each variable~$z_i$ becomes a vertex in~$V$.
We also add an auxiliary
vertex~$s$ to $V$. For each inequality $z_i - z_j \le w_{ij}$ we
add an arc of length~$w_{ij}$ from~$z_j$ to $z_i$.
Finally, for each~$i$ we add an arc~$s \to z_i$ of zero length.

Clearly, the number of arcs in the resulting graph is linear in
the number of constraints of the original system of inequalities. We
invoke Ford--Bellman's shortest-path algorithm
starting from the vertex~$s$. This algorithm runs in~$O(VE)$ time
and either finds a cycle of negative length or computes the
distances from the origin~$s$ to all vertices reachable from
$s$.

Suppose there is a cycle of negative length in $G$. Then it cannot
pass through origin~$s$ since it has no incoming arcs. Hence each
of the arcs of the cycle corresponds to some inequality.
Summing up these inequalities we get a contradiction showing that
the set of inequalities is inconsistent. Otherwise let~$d(u)$ be
the distance from the origin~$s$ to a vertex~$u$. Then triangle
inequality shows that the distances~$d(u)$ obey all the difference
constraints. Moreover, these distances are integers (since the
lengths~$w_{ij}$ are integers).

The total running time of the algorithm is $O(mn (m + n) (1 +
U_{00}))$ (recall that $V = O(m + n)$, $E = O(mn)$ and there are
$1 + U_{00}$ possible values of $A_{00}$). This time bound can be
improved a bit. One may see that instead of $A_{00}$ we may choose
an arbitrary element $A_{ij}$ instead of $A_{00}$ thus proving
Theorem~\ref{upper_2x2}. The running time is polynomial provided
that the elements of $U$ are given in unary notation. An open
question is if there exists an algorithm whose running time is
$\mathrm{poly}(\log U)$.

\label{upper2x3_section}
\section{NP-completeness of the $2 \times 3$ Window Case With Upper Constraints}

In this section we prove that the problem for $2 \times 3$ windows and upper constraints
in unary notation is NP-complete. More precisely, consider the following relation:
$$
    R = \{ \langle A, S, U \rangle \mid 0 \le A \le U, \Sigma_{23}(A) = S \}.
$$
Here $A$, $S$ and $U$ are matrices of any appropriate size.
This relation corresponds to the language~$L(R)$ consisting of pairs
$\langle S, U \rangle$ for which the problem has a solution:
$$
    L(R) = \{ \langle S, U \rangle \mid \exists A \; \langle A, S, U \rangle \in R \}.
$$
It is clear that~$L(R) \in NP$. We present a
Karp reduction from a 3-coloring problem to $L(R)$ thus proving
the NP-completeness of $L(R)$.

It is convenient to consider a slightly more general form of the problem
by imposing \emph{two-sided constraints} on the elements of matrix~$A$:
\begin{equation}
\label{twosided2x3}
    \begin{array}{l}
        L \le A \le U; \\
        \Sigma_{23}(A) = S.
    \end{array}
\end{equation}

Computationally this problem is not harder than the original one.
Indeed, let~$A = L + X$, where~$X \in M_{mn}$. Then constraints~$L
\le A \le U$ become equivalent to $0 \le X \le U - L$. Thus we have
reduced the problem with two-sided
constraints to the problem with upper constraints~$U'$ and
sums~$S'$, where
\begin{align*}
    U' & = U - L; \\
    S' & = S - \Sigma_{23}(L).
\end{align*}
The matrices~$U'$ and $S'$ can be computed in $O(mn)$ time.

Consider a $(m + 1) \times (3n + 2)$ matrix~$A$ of the form:
\begin{center}
    \renewcommand{\arraystretch}{0}
    \begin{tabular}{|c|c||c|c|c||c|c|c||c|c|c||c|}
        \hline
        \strut 0          & 0       & $z_1$ & $p_1$ & $q_1$ & $z_2$ & $p_2$ & $q_2$ & \ldots \\
        \hline
        \rule{0pt}{1.5pt} &         &       &       &       &       &       &       &        \\
        \hline
        \strut $+x_1$     & $+y_1$  &       &       &       &       &       &       &        \\
        \hline
        \strut $-x_2$     & $-y_2$  &       &       &       &       &       &       &        \\
        \hline
        \strut $+x_3$     & $+y_3$  &       &       &       &       &       &       &        \\
        \hline
        \strut $-x_4$     & $-y_4$  &       &       &       &       &       &       &        \\
        \hline
        \strut \vdots    & \vdots   &       &       &       &       &       &       &        \\
        \hline
    \end{tabular}
\end{center}

The properties $A_{00} = A_{01} = 0$ can be ensured by setting $L_{00} = L_{01} = U_{00} = U_{01} = 0$.
We put $S = 0$ and thus all $2 \times 3$ sums of $A$ are zeros.
Then as in the case of $2 \times 2$ windows one may see that for every~$i, j \ge 1$
\begin{equation*}
    \begin{array}{lcl}
        A_{i,3j}   & = & (-1)^{i+1} (x_i - p_j);     \\
        A_{i,3j+1} & = & (-1)^{i+1} (y_i - q_j);     \\
        A_{i,3j+2} & = & (-1)^i (x_i + y_i + z_j).
    \end{array}
\end{equation*}

Therefore, any system of two-sided constraints on the
values
\begin{equation}
\label{constr2x3}
    \begin{array}{l}
        x_i, \: y_j, \: z_j, \: p_j, \: q_j, \\
        x_i - p_j,               \\
        y_i - q_j,               \\
        x_i + y_i + z_j
    \end{array}
\end{equation}
may be reduced to~$L(R)$.

Note that these expressions are of some very special form
(variables are divided into five groups and only some
combinations are allowed). However, it turns out that any system
of two-sided constraints on sums of at most three variables can
be reduced to this special case. 

Using
variables~$p_j$, we can
represent an equation~$x_\alpha = x_\beta$
 (for arbitrary
$\alpha, \beta$) as follows:
\begin{equation*}
    \begin{array}{ccccc}
        0 & \le & x_\alpha - p_j & \le 0; \\
        0 & \le & x_\beta -  p_j  & \le 0.
    \end{array}
\end{equation*}
(we use a ``fresh'' index $j$ for each equation).
Except for that, we do not use variables $p_j$.
The equations~$y_\alpha = y_\beta$ can be expressed in a similar way
using $q_j$.

Now we show how to write an equation $x_\alpha = y_\beta$ for
arbitrary $\alpha$, $\beta$.
Again we choose ``fresh'' indices~$i$, $j$ and $k$ and write
\begin{align*}
    x_\alpha&= x_i; \\
    y_i&= 0; \\
    x_i + y_i + z_j&= 0; \\
    y_\beta&= y_k; \\
    x_k&= 0; \\
    x_k + y_k + z_j&= 0.
\end{align*}

Equation $z_\alpha = z_\beta$ becomes
\begin{align*}
    x_i&= 0; \\
    x_i + y_i + z_\alpha&= 0; \\
    y_j&= 0; \\
    x_j + y_j + z_\beta&= 0; \\
    y_i&= x_j.
\end{align*}
with ``fresh'' indices $i$ and $j$.

The last issue is an equation $x_\alpha = z_\beta$. Consider
``fresh'' indices $i$, $j$ and write
\begin{align*}
    x_i&= x_\alpha; \\
    x_j&= 0; \\
    x_j + y_j + z_\beta&= 0; \\
    y_i&= y_j; \\
    z_0&= 0; \\
    x_i + y_i + z_0&= 0.
\end{align*}
(We may use the same variable~$z_0$ in all such equations.)
Now all variable groups $x_i$, $y_i$ and $z_j$ have become fully symmetric
and a two-sided constraint may be enforced for a sum of arbitrary two or three variables
as required.

Consider an undirected graph~$G = \langle V, E \rangle$. A
\emph{valid 3-coloring} of $G$ assigns one of three colors
to each vertex of $G$ in such a way that no edge
connects the vertices of the same color. The \emph{graph
3-coloring problem} is to find a valid 3-coloring of $G$ or
establish that it does not exist. The corresponding language
$$
    3\mbox{-}COL = \{ G \mid \mbox{ graph~$G$ admits a valid 3-coloring } \}
$$
is known to be NP-complete (see~\cite{CorLeiRiv1999}).

This problem can be stated as an integer program in the following way.
Assign three integer variables~$x_v$, $y_v$, $z_v$ (corresponding to
three possible colors) to each vertex of $G$.
Each of variables should be either 0 or 1:
$$
    0  \le  x_v  \le  1, \: 0  \le  y_v  \le  1, \: 0  \le  z_v  \le  1.
$$
Since each vertex should be assigned a color
$$
    x_v + y_v + z_v = 1.
$$
The requirement that no edge connects the vertices of the same color
produces the following set of inequalities for each edge $uv \in E$:
\begin{eqnarray*}
    x_u + x_v \le 1; \\
    y_u + y_v \le 1; \\
    z_u + z_v \le 1.
\end{eqnarray*}

All these inequalities are constraints on the sum of at most
three variables. Thus these inequalities are equivalent to some $2 \times 3$
problem with two-sided constraints. Clearly this reduction can be
performed in polynomial time and produces matrices $L$, $U$ and $S$
of polynomial size. Thus we have obtained the proof of
Theorem~\ref{upper_2x3}.

\nocite{*}
\bibliographystyle{plain}
\bibliography{main}

\end{document}